\documentclass[10pt,a4paper,reqno]{amsart}
\usepackage{graphicx}

\newtheorem{theorem}{Theorem}
\newtheorem{lemma}{Lemma}

\newcommand{\be}{\begin{equation}}
\newcommand{\ee}{\end{equation}}

\begin{document}
\title[A note on the squeezing function]{A note on the squeezing function} %

\author[A. Yu. Solynin]{Alexander Yu. Solynin}
\address{Department of Mathematics and Statistics, Texas Tech
University, Box 41042, Lubbock, Texas 79409}
\email{alex.solynin@ttu.edu}
\begin{abstract}
The squeezing problem on $\mathbb{C}$ can be stated as follows.
Suppose that $\Omega$ is a multiply connected domain in the unit
disk $\mathbb{D}$ containing the origin $z=0$. How far can the
boundary of $\Omega$  be pushed from the origin by an injective
holomorphic function $f:\Omega\to \mathbb{D}$  keeping the origin
fixed?

In this note, we discuss recent results on this problem obtained
by Ng, Tang and Tsai (Math. Anal. 2020) and by Gumenyuk and Roth
(arXiv:2011.13734, 2020) and also prove few new results using a
method suggested in one of our previous papers (Zapiski Nauchn.
Sem. POMI 1993).

\end{abstract}
\maketitle

\setcounter{equation}{0}  %

\noindent %
\textbf{The squeezing problem.} The squeezing function
$S_\Omega(z)$ of a planar domain $\Omega$ is defined as follows.
Suppose that $\Omega\subset \mathbb{C}$ is such that there is an
injective holomorphic  function $f(z)$ from $\Omega$ to the unit
disk $\mathbb{D}=\{z:\,|z|<1\}$. Let $\mathcal{U}(\Omega)$ be the
class of all such functions. Then $S_\Omega:\Omega\to \mathbb{R}$
is defined by %
\begin{equation}  \label{1} %
S_\Omega(z)=\sup\{{\mbox{dist}}(0,\partial f(\Omega)):\,f\in
\mathcal{U}(\Omega),\,f(z)=0\}. %
\end{equation} %

According to \cite{NTT}, the squeezing function was first
introduced in 2012 by Dong, Guan and Zhang \cite{NTT} but the
concept itself goes back to the work of Liu, Sun and Yau
\cite{LSY}. These authors defined $S_\Omega$ and used it in a more
general setting, namely, for classes of injective holomorphic
mappings defined on domains $\Omega$ in $\mathbb{C}^n$, $n\ge 1$.

The \emph{squeezing problem}, i.e. the problem to find or
characterize $S_\Omega(z)$, is a difficult task, even in complex
dimension~$1$. Any function $f\in \mathcal{U}(\Omega)$ such that
$f(z)=0$ and $S_\Omega(z)={\mbox{dist}}(0,\partial f(\Omega))$
will be called extremal for the squeezing problem (for the point
$z\in \Omega$). In this case, the image $f(\Omega)$ will be called
an extremal domain.

For doubly connected domains $\Omega\subset \mathbb{C}$, the
squeezing function was identified by Ng, Tang and
Tsai in 2020 \cite{NTT}. 
To state the main result of \cite{NTT}, we first introduce
necessary notations. As well known, see \cite[Theorem
5.4]{Jenkins}, every domain $\Omega\subset \mathbb{C}$ with a
finite number of boundary continua $\gamma_1,\ldots,\gamma_n$,
$n\ge 2$, can be mapped by $w=\varphi_{a,k}(z)$, $a\in \Omega$,
$1\le k\le n$, conformally on $\mathbb{D}$ slit along arcs
$C_{j,k}=C_{j,k}(a)$, $j\not= k$, on the circles of radii
$0<r_{j,k}=r_{j,k}(a)<1$ centered at $0$ so that
$\varphi_{a,k}(a)=0$, a non-degenerate boundary continuum
$\gamma_k$ corresponds to the unit circle $\mathbb{T}=\partial
\mathbb{D}$ and, for $j\not= k$, $\gamma_j$ corresponds to
$C_{j,k}$. Under the additional normalization
$\varphi'_{a,k}(a)>0$, that we assume in what follows,  the
mapping function $\varphi_{a,k}$ is uniquely determined. If
$\gamma_k$ is a singleton, we put $r_{j,k}(a)=0$ for all $a \in
\Omega$ and all $j\not=k$. Without the uniqueness statement, the
mapping property, discussed above, remains true for domains of any
connectivity.

\begin{theorem}[\cite{NTT},\cite{GR}] \label{Theorem 1}  %
Let $\Omega\subset \mathbb{C}$ be a doubly connected domain with
boundary components $\gamma_1$ and $\gamma_2$, at least one of
which is non-degenerate.
Then %
\begin{equation} \label{2}%
S_\Omega(z)=\max\{r_{1,2}(z),r_{2,1}(z)\}.
\end{equation}  %
If $r_{1,2}(z)<r_{2,1}(z)$, then $\varphi_{z,1}$ is the unique (up
to rotation about the origin) extremal function for the squeezing
problem; if $r_{1,2}(z)>r_{2,1}(z)$, then $\varphi_{z,2}$ is the
unique  extremal function; and if $r_{1,2}(z)=r_{2,1}(z)$, then
both $\varphi_{z,1}$ and $\varphi_{z,2}$
are extremal. %
\end{theorem}  %

This theorem, except for the uniqueness part, was proved in
\cite{NTT} using the Loewner differential equation and tricky
calculations
with special functions. A simpler proof 
based on the potential theory, which also includes the proof of
the uniqueness statement, was presented in \cite{GR}.

In this note, we first show that Theorem~1 is immediate from
Theorem~1.2 in our 1993 paper \cite{S1}. Then we discuss how the
method used in \cite{S1} can be applied to study the squeezing
problem in a more general setting.

Let $0<r_1\le r_2<1$ and let $E$ be a compact subset in the
closure $\bar{A}(r_1,r_2)$ of the annulus
$A(r_1,r_2)=\{z:\,r_1<|z|<r_2\}$. Let
$\mathcal{D}(E)=\{(D_1,D_2)\}$ be the set of pairs $(D_1,D_2)$ of
non-overlapping domains in $\mathbb{D}\setminus E$, where $D_1$ is
a simply connected domain containing the origin and $D_2$ is a
doubly connected domain separating $E$ from $\mathbb{T}$. In what
follows,   $m(D, z_0)$ stands for the reduced module of a simply
connected domain $D$ with respect to the point $z_0\in D$ and
$m(D)$ stands for the module of a doubly connected domain $D$. For
the definitions and properties of these moduli, we refer to
Jenkins's monograph \cite{Jenkins} as the primary source and also
to \cite{Kuzmina}, \cite{Dubinin}, and \cite{S2} . Figure~1
illustrates our notations introduced above.

Consider the following \emph{module problem}.

\noindent %
\textbf{Problem M.} \emph{Given nonnegative numbers $\alpha_1$ and
$\alpha_2$, at least one of which is positive, identify all pairs
$(D_1^*,D_2^*)\in \mathcal{D}(E)$, which maximize the weighted sum
of moduli
\begin{equation}  \label{3} %
\alpha_1^2 m(D_1,0)+\alpha_2^2 m(D_2)
\end{equation}  %
over the set $\mathcal{D}(E)$.}

\begin{theorem}[{\mbox{\cite[Theorem 1.2]{S1}}}] \label{Theorem 2}  %
(1) There is a unique pair $(D_1^*,D_2^*)\in \mathcal{D}(E)$
maximizing the sum (\ref{3}) over the class $\mathcal{D}(E)$. The
domains $D_1^*$ and $D_2^*$ are, respectively, a circle domain and
a ring domain of a quadratic differential $Q(z)\,dz^2$ defined on
$\mathbb{D}\setminus E$, which is positive on $\mathbb{T}$ and has
a second order pole with circular local structure of trajectories
at $z=0$.

(2) Let $L(E)$ denote the \emph{free boundary} of the module
problem; i.e. $L(E)=(\partial D_1^*\cup \partial D_2^*)\cap
(\mathbb{D}\setminus E)$. Then $L(E)$ consists of arcs of critical
trajectories of $Q(z)\,dz^2$ and their endpoints in
$\mathbb{D}\setminus E$.

(3) The following holds: (a) if $\alpha_1\le \alpha_2$, then
$L(E)\subset \bar{\mathbb{D}}_{r_2}$, where
$\mathbb{D}_r=\{z:\,|z|<r\}$; (b) if $\alpha_1\ge \alpha_2$, then
$L(E)\subset \bar{A}(r_1,1)$; (c) if
$\alpha_1=\alpha_2$, then $L(E)\subset \bar{A}(r_1,r_2)$. %
\end{theorem}  %

We note that Problem~M is a particular case of Jenkins's problem
on extremal partitioning of Riemann surfaces, see \cite{Jenkins2},
\cite{Kuzmina}, \cite{S2}. Therefore, parts (1) and (2) of this
theorem follow from Jenkins's theorem on extremal partitioning;
see \cite[Theorem~1]{Jenkins2}. Part (3), which is essential for
this paper, was proved in \cite{S1}. For convenience of the
readers and because this result seems useful in the study of the
squeezing function~(\ref{1}), a version of this proof will be
presented at the end of this note.

\medskip
\begin{figure} \label{ModuleProblem1}%
\hspace{-2.5cm} %
\begin{minipage}{1.0\textwidth}%
$$\includegraphics[scale=0.5,angle=0]{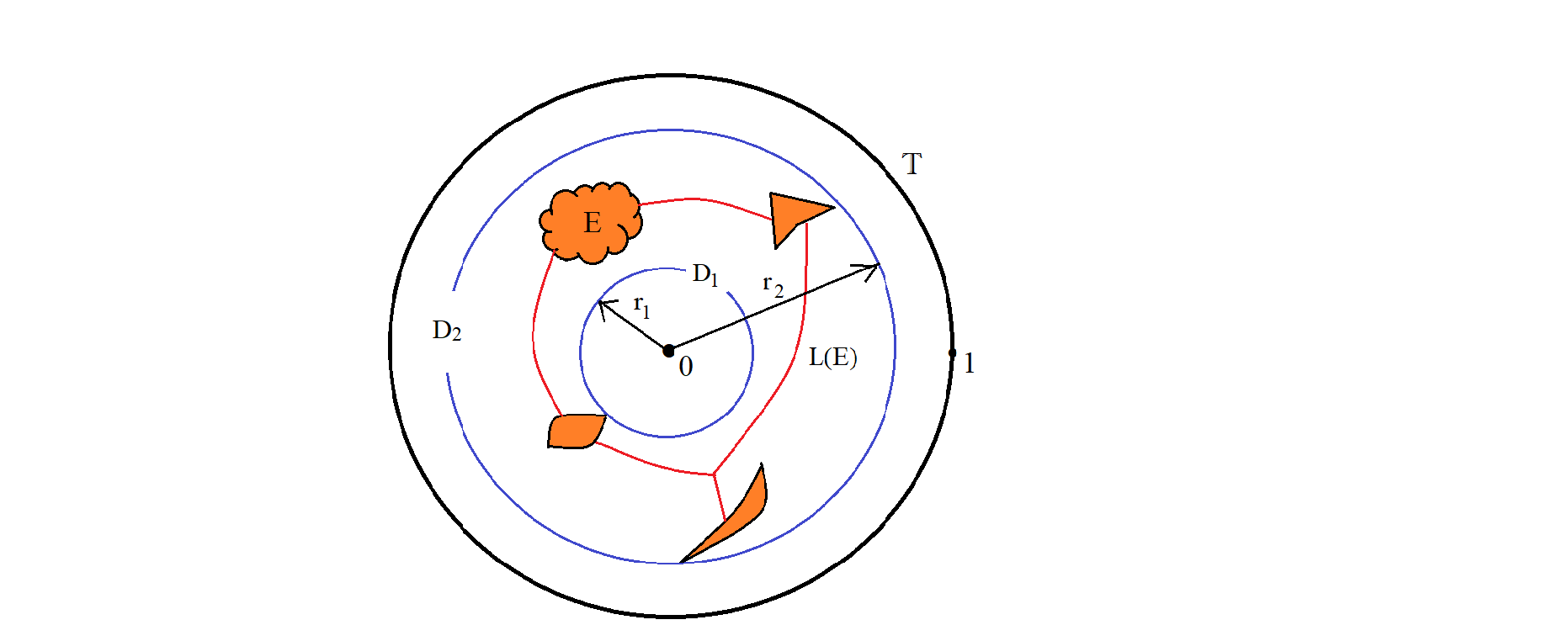}
$$
\vspace{-1cm} %
\end{minipage}
\caption{Extremal partitioning into two domains.}%
\end{figure}

\medskip

\noindent %
\emph{Proof of Theorem~\ref{Theorem 1}.} Using an auxiliary
conformal mapping, we may assume without loss of generality that
$\Omega$ is the unit disk $\mathbb{D}$ slit along a proper arc
$C_{1,2}$ of the circle $\{z:\,|z|=r_{1,2}\}$, that $z=0$, and
that $r_{1,2}\ge r_{2,1}$. Consider Problem~M with
$\alpha_1=\alpha_2=1$ for the set $E=C_{1,2}$. We call it
Problem~$P_1$. It follows from Gr\"{o}tzsch's lemma
\cite[Theorem~2.6]{Jenkins} (or from Theorem~2(3) above) that
$(\mathbb{D}_{r_{1,2}},A(r_{1,2},1))$ is the extremal pair of
domains for Problem~$P_1$ and the complementary arc
$C'_{1,2}=\{z:\,|z|=r_{1,2}\}\setminus C_{1,2}$ is the free
boundary of Problem~$P_1$.

Now, suppose, by contradiction, that $S_\Omega(0)>r_{1,2}$.  The
latter implies that there is a function $f\in
\mathcal{U}(\mathbb{D}\setminus C_{1,2})$ such that $f(0)=0$ and
${\mbox{dist}}(0,f(C_{1,2}))>r_{1,2}$. Let
$\rho_1=\inf\{|z|:\,z\in f(C_{1,2})\}$, $\rho_2=\sup\{|z|:\,z\in
f(C_{1,2})\}$. Consider Problem~M, again with
$\alpha_1=\alpha_2=1$, for the set $f(C_{1,2})\subset
\bar{A}(\rho_1,\rho_2)\}$. We call it Problem~$P_2$. Since
Jenkins's problem on extremal partitioning is conformally
invariant (in the sense that conformal mappings preserve extremal
configurations), it follows that the pair
$(f(\mathbb{D}_{r_{1,2}}),f(A(r_{1,2},1)))$ is extremal for
Problem~$P_2$. Also, conformal mappings preserve the module of a
doubly connected
domain. Hence, %
\begin{equation} \label{4} %
m(f(A(r_{1,2},1)))=m(A(r_{1,2},1))=-\frac{1}{2\pi}\log r_{1,2}.
\end{equation} %
 Now, it follows from Theorem~2(3) that the free
boundary of Problem~$P_2$ lies in the annulus
$\bar{A}(\rho_1,\rho_2)$. This implies that
$f(A(r_{1,2},1))\subset A(\rho_1,1)$. Since the module of a doubly
connected domain  increases under expansion, it follows that
$m(f(A(r_{1,2},1)))\le m(A(\rho_1,1))=-\frac{1}{2\pi}\log \rho_1$.
Since $r_{1,2}<\rho_1$, the latter inequality contradicts
(\ref{4}). Thus, our assumption leads to a contradiction and
therefore we must have $S_\Omega(0)=r_{1,2}$.

It remains to show that rotations about the origin are the only
extremal functions for the problem under consideration. If $f$ is
not a rotation, then our argument above implies that
$f(A(r_{1,2},1))$ is a proper subdomain of $A(r_{1,2},1)$.
Therefore, $m(f(A(r_{1,2},1)))<m(A(r_{1,2},1))=-\frac{1}{2\pi}\log
r_{1,2}$, which again contradicts equation (\ref{4}). The proof is
 complete. %
  \hfill $\Box$

 \bigskip

 \noindent %
\textbf{ More results and questions.} Let $\Omega$ be a domain (of
any connectivity) and let $\gamma$ be a nondegenerate boundary
continuum of $\Omega$ separated from the rest of $\partial
\Omega$, i.e. such that $\gamma\cap
(\overline{\partial\Omega\setminus \gamma})=\emptyset$. If
$\Omega$ is finitely connected, the latter separation property
always holds. Then, for each $z\in \Omega$, we can consider a
non-empty class $\mathcal{U}_\gamma(\Omega)$ of functions $f\in
\mathcal{U}(\Omega)$ such that $f(\gamma)=\mathbb{T}$, in the
sense of boundary correspondence. Then the function
$S_{\Omega,\gamma}:\Omega\to \mathbb{R}$ defined by
\begin{equation}  \label{5} %
S_{\Omega,\gamma}(z)=\sup\{{\mbox{dist}}(0,\partial
f(\Omega)):\,f\in
\mathcal{U}\gamma(\Omega),\,f(z)=0\} %
\end{equation} %
can be thought as the squeezing function toward the boundary
continuum $\gamma$. For this function we have the following
result.

\begin{lemma} \label{Lemma 1}  %
Let $\Omega$ and $\gamma$ be as above and let $a\in \Omega$.
Suppose that there is a function $f_{a,\gamma}\in
\mathcal{U}_\gamma(\Omega)$ such that $f_{a,\gamma}(a)=0$ and
$\mathbb{D}\setminus f_{a,\gamma}(\Omega)\subset \{z:\,|z|=r\}$,
$0<r<1$. Then $S_{\Omega,\gamma}(a)=r$.

Furthermore, $f_{a,\gamma}$ is a unique (up to rotation about $0$)
function in $\mathcal{U}_\gamma(\Omega)$ that is extremal for the
problem on the squeezing toward $\gamma$ for the point $a$.%
\end{lemma}  %

\noindent %
\emph{Proof.} Since the squeezing problem is conformally
invariant, we may assume that $\Omega$ is the disk $D$ slit along
a proper compact subset of the circle $T_r=\{z:\,|z|=r\}$,
$0<r<1$. Then our argument used in the proof of Theorem~3 shows
that if $f\in \mathcal{U}_\mathbb{T}(\Omega)$ such that $f(0)=0$
is extremal for the squeezing toward $\mathbb{T}$ problem, then
$\mathbb{D}\setminus f(\Omega)\subset \mathbb{T}_r$. Furthermore,
same argument shows that $f(\mathbb{D}_r)=\mathbb{D}_r$. Since
$f(0)=0$, the latter implies that $f$ is a rotation about the
origin. %
\hfill $\Box$

\medskip

Below, we assume that $\Omega$ is a finitely-connected domain with
$n\ge 3$ nondegenerate boundary continua
$\gamma_1,\ldots,\gamma_n$. In this case, the following properties
are either exist in the literature or easy to prove:
\begin{enumerate}
\item[(1)] \textbf{Existence of extremal functions.} For each
$z\in \Omega$ and $1\le k\le n$, there is an extremal function for
the problem on squeezing toward $\gamma_k$; i.e. a function
$f_{z,k}\in \mathcal{U}_{\gamma_k}(\Omega)$ such that
$f_{z,k}(z)=0$ and
$S_{\Omega,\gamma_k}(z)={\mbox{dist}}(0,\partial
f_{z,k}(\Omega))$. Therefore, for each $z\in \Omega$, there is a
function $f_\Omega\in \mathcal{U}(\Omega)$
extremal for the squeezing problem on $\Omega$. %
\item[(2)] \textbf{Continuity.} The functions $S_\Omega(z)$ and
$S_{\Omega,\gamma_k}(z)$, $k=1,\ldots,n$, are continuous on
$\Omega$. %
\item[(3)] \textbf{Monotonicity of $S_{\Omega,\gamma_n}(z)$ with
respect to $\gamma_n$. } For a fixed $z\in \Omega$,
$S_{\Omega,\gamma_n}(z)$ is monotone with respect to $\gamma_n$ in
the following sense. Let $\Omega'\not= \Omega$ be a finitely
connected domain with boundary continua
$\gamma'_1,\ldots,\gamma'_{n'}$, $2\le n'\le n$, such that $a\in
\Omega' \subset \Omega$ and
$\gamma'_k=\gamma_k$ for $k=1,\ldots,n'-1$. Then %
$$
S_{\Omega',\gamma'_{n'}}(a)>S_{\Omega,\gamma_n}(a).
$$  %
\item[(4)] \textbf{Non-monotonicity of $S_\Omega(z)$ with respect
to $\Omega$.} For a fixed $z\in \Omega$, $S_\Omega(z)$ is not
monotone as a function
of the domain. %
\item[(5)] \textbf{Boundary values of $S_\Omega(z)$.}
$S_{\Omega,\gamma_k}(z)\to 1$ as $z\to \gamma_k$ and therefore
$S_\Omega(z)\to 1$ as $z\to
\partial \Omega$.
\item[(6)] \textbf{Boundary values of $S_{\Omega,\gamma_k}(z)$ on
$\gamma_j$, $j\not= k$.} For each $k$ and $j\not= k$, there is
$0<c_{j,k}<1$ such that %
$$
\limsup_{z\to \gamma_j} S_{\Omega,\gamma_k}(z)\le c_{j,k}.
$$
\item[(7)] \textbf{Separation property 1.} If $c$, $0<c<1$, is
sufficiently close to $1$, then the level set $\{z\in
\Omega:\,S_{\Omega,\gamma_k}(z)=c\}$ separates $\gamma_k$ from
$\partial \Omega\setminus \gamma_k$. %
\item[(8)] \textbf{Separation property 2.} Let $j\not=k$. The set
$\{z\in \Omega:\,S_{\Omega,\gamma_j}(z)=S_{\Omega,\gamma_k}(z)\}$
separates $\gamma_j$ from $\gamma_k$ inside $\Omega$. %
\end{enumerate} %

\medskip

(1) The existence of extremal functions and therefore existence of
extremal domains as well was established in
\cite[Theorem~2.1]{DGZ}.

(2) To prove the continuity property, one can use extremal
functions $f_{z,\Omega}(\zeta)$ composed with the Moebius
automorphisms of the unit disk $\mathbb{D}$. The details are left
to the interested reader.

(3) To prove the monotonicity property of the squeezing function
$S_{\Omega,\gamma_n}(z)$ defined by (\ref{5}), we consider the
composition $\varphi\circ f_{\Omega,\gamma_n}$ of the extremal
function $f_{\Omega,\gamma_n}$ with the function $\varphi$ that is
the Riemann mapping function from a simply connected domain with
the boundary continuum $\gamma'_{n'}$ such that $\varphi(0)=0$.
Then the desired result follows from the Schwarz's lemma.

(4) To show that the monotonicity with respect to $\Omega$ is
absent, consider two simple examples. Let $\Omega$ be a circularly
slit disk $\mathbb{D}\setminus \{re^{i\theta}:\, |\theta|\le
\alpha\}$ with $0<r<1$, $0<\alpha<\pi$. For $0<r_1<r$, $r<r_2<1$,
let $\Omega_1=\Omega\setminus [r_1,r]$, $\Omega_2=\Omega\setminus
[r,r_2]$. Then the argument involving Theorem~2, as it was used in
the proof of Theorem~1, shows that
$S_{\Omega_1}(0)<S_\Omega(0)<S_{\Omega_2}(0)$. Therefore, the
monotonicity property of $S_\Omega(z)$ as a function of the domain
does not hold, in general.

(5) To show that $S_{\Omega,\gamma_k}(z)\to 1$ as $z\to \gamma_k$,
we may assume that $\Omega\subset \mathbb{D}$ and
$\gamma_k=\mathbb{T}$. Then the Moebius mapping
$\varphi(\zeta)=(\zeta-z)/(1-\bar{z}\zeta)$ is in
$\mathcal{U}_{\gamma_k}(\Omega)$ and it is easy to see that that
${\mbox{dist}}(0,\partial \varphi(\Omega))\to 1$ as $|z|\to 1$.

(6) Let $l$ be a Jordan curve in $\Omega$ separating   $\gamma_j$
from $\partial \Omega\setminus \gamma_j$ and let $\Omega_l$ be a
doubly connected domain with boundary components $\gamma_j$ and
$l$. It follows from the property (3) that
$S_\Omega(z)<S_{\Omega_l}(z)$ for all $z\in\Omega_l$.

Fix $z_0\in \Omega_l$. Let $m=m(D_l)$ be the module of the doubly
connected domain $D_l$. Then there is a function $\varphi$ that
maps $\Omega_l$ conformally onto the annulus $A(s,1)$ with
$s=e^{-2\pi m}$. In addition, we may assume that $\varphi(z_0)=a$,
$s<a<1$. Furthermore, there exists a function $\psi$, such that
$\psi(a)=0$, which maps $A(s,1)$ conformally onto the disk
$\mathbb{D}$ slit along and arc $C$ on a circle $\mathbb{T}_\rho$,
with $0<\rho<1$, where $\rho=\rho(a)$ depends on $a$. The mapping
$\psi$ is one of the well-studied canonical mappings. A particular
property of $\psi$, we need here, was proved by E.~Reich and
S.~E.~Warschawski in 1960 \cite[Lemma 3]{RW}. These authors showed
that $\rho(a)=a$. This result implies that $a\to s$ as $z_0\to
\gamma_j$. It follows from Theorem~1 that
$S_{\Omega_l}(z_0)=\rho(a)=a$. This, being combined with
property (3), implies that %
$$
\limsup_{z_0\to \gamma_j}\le \lim_{a\to s}a=s=e^{-2\pi m}<1,
$$
as required.

(7) Let $c_k=\max_{j\not = k}\{c_{j,k}\}$, where $c_{j,k}$
introduced in part (6) above. Take $c$, $c_k<c<1$, and consider
the level set $L(c)=\{z\in \Omega: S_{\Omega,\gamma_k}(z)=c\}$.
Since $S_{\Omega,\gamma_k}$ is continuous on $\Omega$ and has
boundary value $1$ on $\gamma_k$ and boundary values less than $c$
on other boundary continua, it follows that $L(c)$ separates
$\gamma_k$ from the rest of $\partial \Omega$.

(8) The same argument, which was used in part (7), can be used to
prove the separation property in question as well.

\medskip

Let $\Omega$ be an $n$-connected domain as above. Then, for each
$z\in \Omega$, there are $n$ functions $\varphi_{z,k}$, $1\le k\le
n$, and for each $k$ there are $n-1$ radii $r_{j,k}=r_{j,k}(z)$ of
circular arcs $C_{j,k}=C_{j,k}(z)$, as described above. An example
of the domain $\Omega$ and its image under one of the mappings
$\varphi_{z,k}$ are shown in Figure~2. It was conjectured in
\cite{NTT} that %

\begin{equation}  \label{6} %
S_\Omega(z)=\max_k \min_{j\not=k} r_{j,k}(z). %
\end{equation} %

\medskip
\begin{figure} \label{SqueezingProblem1}%
\hspace{-2.5cm} %
\begin{minipage}{1.0\textwidth}%
$$\includegraphics[scale=0.5,angle=0]{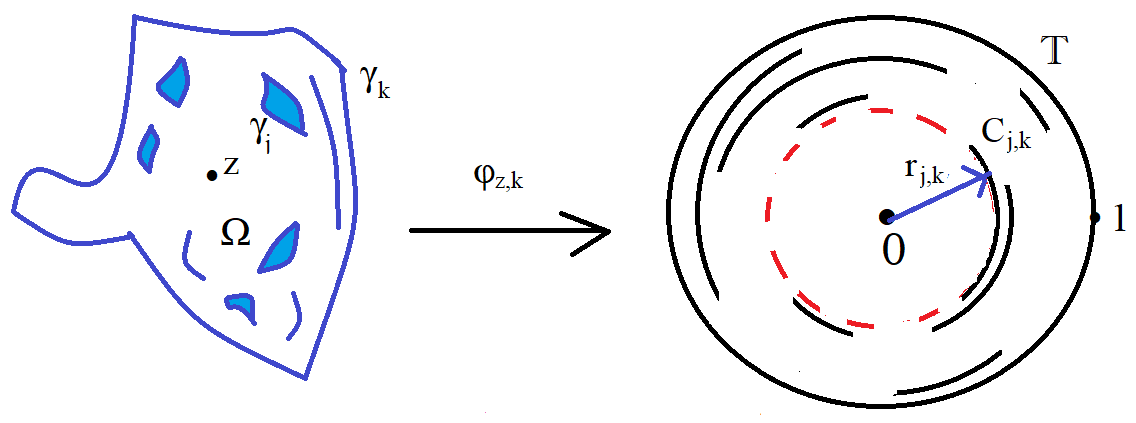}
$$
\vspace{-0.5cm} %
\end{minipage}
\caption{Mapping onto a circularly slit disk.}%
\end{figure}

\medskip

This conjecture, though it sounds plausible,  was quickly
disproved by Gumenyuk and Roth \cite{GR}, whose engineering
construction shows that for every $n\ge 3$ there is an
$n$-connected domain $\Omega$, that is not a circularly slit disk,
and there is a point $z\in \Omega$ such that $S_\Omega(z)$ is
strictly greater than the right-hand side of (\ref{6}). This shows
that, in general, the family of domains extremal for the squeezing
problem is not limited to the set of circularly slit disks. In
\cite{NTT} and \cite{GR}, the question whether or not there exist
circularly slit disks of connectivity $n\ge 3$ that are extremal
for the squeezing problem was left open. Thus, we complement
results in \cite{NTT} and \cite{GR} with the following theorem.

\begin{theorem} \label{Theorem 3}  %
For each $n\ge 3$, there is a domain $\Omega$, that is a
circularly slit disk of connectivity $n$, such that $\Omega$ is
the extremal domain, unique up to rotation about the origin, for
the squeezing problem for the point $z=0$.
\end{theorem}  %
\noindent %
\emph{Proof.} Let $n\ge 3$, $0<r<1$, and $0<\alpha<\pi/(n-1)$. Let
$\gamma_k=\{z=re^{i\theta}:\,|\theta-2\pi(k-1)/(n-1)|\le \alpha\}$
for $1\le k\le n-1$ and let $\gamma_n=\mathbb{T}$. Let
$\Omega=\Omega(n,r,\alpha)$ denote the disk $\mathbb{D}$ slit
along the arcs $\Gamma_k$, $k=1,\ldots,n-1$.

It follows from Lemma~1 that $S_{\Omega,\gamma_n}(0)=r$.
Furthermore, since $\Omega$ possesses $(n-1)$-fold rotational
symmetry, it follows that there is $\rho$, $0<\rho<1$, such that %
$$ 
S_{\Omega,\gamma_k}(0)=\rho \quad \quad {\mbox{for
$k=1,\ldots,n-1$.}}
$$ 

We claim that, for fixed $n$ and $r$, $\rho=\rho(\alpha)\to 0$ as
$\alpha\to 0$. To emphasize dependence on $\alpha$, we use
notations $\Omega=\Omega(\alpha)$, $\gamma_k=\gamma_k(\alpha)$,
etc. To prove the claim, we consider a circle
$\Gamma_\varepsilon=\{z:\,|z-r|=\varepsilon\}$ with
$\varepsilon>0$ sufficiently small. If $\alpha$ is small, then
$\Gamma_\varepsilon$ separates $\gamma_1(\alpha)$  from all other
boundary components of $\Omega(\alpha)$. Let
$D_\varepsilon(\alpha)$ denote the doubly connected domain bounded
by $\Gamma_\varepsilon$ and $\gamma_1(\alpha)$. It is clear that
$m(D_\varepsilon(\alpha))\to \infty$ as $\alpha\to 0$. Since the
module of a doubly connected domain is conformally invariant, it
follows that for every $\delta>0$ there is $\alpha_0>0$ such that
$f(\Gamma_\varepsilon)\subset \mathbb{D}_\delta$, whenever
$0<\alpha<\alpha_0$ and $f\in
\mathcal{U}_{\gamma_1(\alpha)}(\Omega(\alpha))$ is such that
$f(0)=0$. Since the images $f(\mathbb{T})$ and
$f(\gamma_k(\alpha))$, $k=2,\ldots,n-1$, lie in the bounded
component of $\mathbb{C}\setminus f(\Gamma_\varepsilon)$, the
claim follows. Therefore, $\rho(\alpha)<r$ if $\alpha$ is small
enough. Hence $S_{\Omega(\alpha)}(0)=r$ for all such $\alpha$ and
$\Omega(\alpha)$ is the only extremal domain up to rotation about
the origin.
The proof is complete. 
\hfill $\Box$%

\medskip

\noindent %
\textbf{Example.} Given $a>1$, let
$\gamma_k=\gamma_k(a)=\{z=te^{2\pi i(k-1)/3}:\,1\le t\le a\}$,
$k=1,2,3$. Then $\Omega=\overline{\mathbb{C}}\setminus
\cup_{k=1}^3 \gamma_k$ is a domain on $\overline{\mathbb{C}}$ of
connectivity $3$. Let $\varphi_{0,1}$ map $\Omega$ conformally
onto a circularly slit disk $\mathbb{D}$ such that
$\varphi_{0,1}(\gamma_1)=\mathbb{T}$,  $\varphi_{0,1}(0)=0$, and
$\varphi'_{0,1}(0)>0$. Since $\Omega$ is invariant under rotation
by angle $2\pi/3$ about $0$, the symmetry principle implies that
$\varphi_{0,1}(\gamma_2)$ and $\varphi_{0,1}(\gamma_3)$ are
circular arcs lying on the same circle symmetrically to each other
with respect to the real axis. Therefore, it follows from
Theorem~3 that $\varphi_{0,1}$ is extremal for the squeezing
toward $\gamma_1$ problem for the point $0$. Since $\Omega$ is
invariant under rotations by the angle $2\pi/3$, it follows that
each of the other two functions, $\varphi_{0,2}$ and
$\varphi_{0,3}$, where $\varphi_{0,k}$ maps $\gamma_k$ to the unit
circle $\mathbb{T}$ such that $\varphi_{0,k}(0)=0$, $k=2,3$, is
also extremal for the squeezing problem toward, respectively,
$\gamma_2$ or $\gamma_3$, for the point $0$. Moreover,
$S_{\Omega,\gamma_1}(0)=S_{\Omega,\gamma_2}(0)=S_{\Omega,\gamma_3}(0)$.

Thus, this example shows that there are domains $\Omega$ of
connectivity $n> 2$ and points $z_0\in \Omega$, such that every
function $f$ extremal for the squeezing problem for $\Omega$ and
$z_0$ maps $\Omega$ onto a circularly slit disk.

\medskip

\noindent %
\textbf{Questions.} We finish this part with six questions for future study. %
\begin{enumerate} %
\item[(1)] As Theorem~1 shows, if $\Omega$ is doubly connected,
then every function $f$ extremal for the squeezing problem for
some point $z_0\in \Omega$ maps $\Omega$ onto a circularly slit
disk. We suspect that this property characterizes doubly connected
domains. So, the question: Is it true that if there is an open
subset $G$ of $\Omega$, such that any function $f$ extremal for
the squeezing problem for some $z\in G$ maps $\Omega$ onto a
circularly slit disk, then $\Omega$ is doubly connected?  %
\item[(2)] A point $z_0$ in a domain $\Omega$ with boundary
continua $\gamma_1,\ldots,\gamma_n$, $n\ge 2$,  is an equilibrium
point for the squeezing problem if
$S_{\Omega,\gamma_1}(z_0)=S_{\Omega,\gamma_2}(z_0)=\ldots=S_{\Omega,\gamma_n}(z_0)$.
An annulus $A(r,1)$ has the circle $\{z:\,|z|=\sqrt{r}\}$ as its
set of equilibrium points. A domain
$D_n(a)=\overline{\mathbb{C}}\setminus \cup_{k=1}^n \{te^{2\pi
i(k-1)/n}:\,1\le t\le a\}$ with $a>1$ and $n\ge 3$ has two
equilibrium points, $z_1=0$ and $z_2=\infty$. Let
$D_n(a,R)=D_n(a)\cap \mathbb{D}_R$, where $R>a$. There is a unique
value $R_0>a$ such that $D_n(a,R)$ has exactly one equilibrium
point $z_0=0$, while, for $R\not= R_0$, $D_n(a,R)$ does not have
equilibrium points.

Is it true that every domain $\Omega$ of connectivity $n\ge 3$ has
at most two equilibrium points? %
\item[(3)] Are there circularly slit disks with boundary continua
$\gamma_n=\mathbb{T}$ and $\gamma_k\subset \mathbb{D}$,
$k=1,\ldots,n$, $n\ge 3$, that are extremal for the squeezing
problem for $z=0$, such that $\gamma_j$ and $\gamma_k$ lie on
different
circles when $j\not= k$? %
\item[(4)]  Are there domains $\Omega\subset \mathbb{D}$ of
connectivity $n+m+1$ with boundary continua
$\gamma_{n+m+1}=\mathbb{T}$ and $\gamma_k\subset \mathbb{D}$,
$k=1,\ldots,n+m$, such that $\gamma_k$ is an arc of a circle
centered at $0$ when $1\le k\le n$, $\gamma_k$ is not a circular
arc when $n+1\le k\le n+m$, that are extremal for the squeezing
problem for $z=0$? %
\item[(5)] Is it possible to characterize domains $\Omega$ of
connectivity $n\ge 3$, that are extremal for the squeezing problem
for a point $z_0\in \Omega$,  in terms of quadratic differentials
similar to characterization of  configurations of domains extremal
for Jenkins's problem \cite{Jenkins2} on extremal partitioning? %
\item[(6)] Our proof of Theorem~2 on the extremal partitioning of
the disk $\mathbb{D}$ with a compact barrier $E\subset \mathbb{D}$
can be extended to higher dimensions for the unit ball
$\mathbb{B}\subset \mathbb{R}^n$, $n\ge 3$, with a compact set
$E\subset \mathbb{B}$.

It would be interesting to know whether such a generalization can
be used to study the squeezing problem in $\mathbb{C}^n$.
\end{enumerate} %

\bigskip

 \noindent %
\textbf{ Proof of Theorem~2.} Problem~M is a particular case of
Jenkins's problem on extremal partitioning \cite{Jenkins2}.
Therefore, parts (1) and (2) of Theorem~2 follow from
\cite[Theorem 1]{Jenkins2}. Furthermore, same Jenkins's theorem
implies that the metric
$\rho_1(z)|dz|=\alpha_1^{-1}|Q(z)|^{1/2}\,|dz|$ is extremal for
the module problem for the family of closed curves $\gamma\subset
D_1^*$ separating $z=0$ from $\partial D_1^*$ and the metric
$\rho_2(z)|dz|=\alpha_2^{-1}|Q(z)|^{1/2}\,|dz|$ is extremal for
the module problem for the family of closed curves $\gamma\subset
D_2^*$ separating $\mathbb{T}$ from $\partial D_2^*\setminus
\mathbb{T}$. Moreover, it implies that %
\begin{equation} \label{8} %
\int_\gamma |Q(z)|^{1,2}\, |dz|\ge \alpha_2 %
\end{equation}  %
for every Jordan curve $\gamma\subset \mathbb{D}\setminus E$,
which separates $\mathbb{T}$ from $E$ and $0$ and also for
$\gamma=\mathbb{T}$.

Let us we consider the case $\alpha_1\le \alpha_2$, assuming that
$E$ consists of a finite number of connected components. In this
case $D_2^*\not= \emptyset$. Otherwise, $\mathbb{T}$ would be on
the boundary of $D_1^*$, which implies that %
$$ %
\int_\mathbb{T} |Q(z)|^{1/2}\,|dz|<\alpha_1\le \alpha_2,
$$ %
contradicting equation (\ref{8}). Furthermore, since $E$ consists
of a finite number of components, it follows that the free
boundary $L(E)$ consists of a finite number of critical
trajectories of $Q(z)\,dz^2$ and their endpoints in
$\mathbb{D}\setminus E$.

Suppose by contradiction that %
\begin{equation}  \label{9}
r_2<\rho=\max \{|z|:\, z\in L(E)\}<1.
\end{equation} %
Since $L(E)$ consists of a finite number of analytic arcs, it
follows that the intersection $L(E)\cap \mathbb{T}_\rho$ consists
of a finite number of points $z_k$, $k=1,\ldots,n$. It follows
from the local structure of trajectories of $Q(z)\,dz^2$ near
critical points, that each of the points $z_k$ is regular. This
implies that, for all sufficiently small  $\varepsilon>0$, the set
$L(E)\cap \bar{A}(\rho-\varepsilon,\rho)$ consists of $n$ disjoint
analytic arcs $s_k$, $k=1,\ldots,n$, such that $s_k$ has its
endpoints on the circle $\mathbb{T}_{\rho-\varepsilon}$ and the
point $z_k$ is an interior point of $s_k$. Let $\tilde{s}_k$
denote the arc symmetric to $s_k$ with respect to the circle
$\mathbb{T}_{\rho-\varepsilon}$. For every $\varepsilon$ small
enough and each $k$, the intersection $\tilde{s}_k\cap L(E)$ is
empty. In this case, there are $n$ simply connected domains
$\Delta_k$, symmetric with respect to the circle
$\mathbb{T}_{\rho-\varepsilon}$, such that $\partial
\Delta_k=s_k\cup \tilde{s}_k$, $k=1,\ldots,n$.

We claim that there is no $k$ such that $\Delta_k\subset D_2^*$.
To prove this claim, suppose that $\Delta_{k_j}\subset D_2^*$ for
$j=1,\ldots, n_1$, $n_1 \le n$. Let $\partial^0 D_2^*=\mathbb{T}$
and $\partial^1 D_2^*$ denote the boundary components of $D_2^*$.
Let $D_2^p$ denote polarization of the doubly connected domain
$D_2^*$ (considered as a condenser whose plates are connected
components of $\mathbb{C}\setminus D_2^*$) with respect to the
circle $\mathbb{T}_{\rho-\varepsilon}$. In the case under
consideration, the polarized domain $D_2^p$ is a doubly connected
domain having $\mathbb{T}$ as one of its boundary components while
the other boundary component of $D_2^p$ is obtained from
$\partial^1D_2^*$ by replacing the arcs  $s_{k_j}$ with the arcs
$\tilde{s}_{k_j}$, $j=1,\ldots, n_1$. For the definition and
properties of polarization, we refer to \cite[Chapter 3]{Dubinin}.
As well known, polarization increases the module of a
doubly connected domain. Thus, in our case, we have %
$$ %
m(D_2^*)<m(D_2^p),
$$ %
with the sign of strict inequality because $D_2^p$ does not
coincide with $D_2^*$ up to reflection with respect to
$\mathbb{T}_{\rho-\varepsilon}$. Therefore,
$$ 
\alpha_1^2m(D_1^*,0)+\alpha_2^2
m(D_2^*)<\alpha_1^2m(D_1^*,0)+\alpha_2^2 m(D_2^p). %
$$ 
Since the pair $(D_1^*,D_2^p)$ is admissible for Problem~M, i.e.
$(D_1^*,D_2^p)\in \mathcal{D}(E)$, inequality (\ref{9})
contradicts our assumption that the pair of domains
$(D_1^*,D_2^*)$ is extremal for Problem~M.

\medskip

Now we consider the case when $\Delta_k\subset D_1^*$ for all
$k=1,\ldots,n$. In this case, we consider a pair $(D_1,D_2)$,
where $D_1=D_1^*\setminus \bar{A}(\rho-\varepsilon,\rho)$ is a
simply connected domain and $D_2=D_2^*\cup
A(\rho-\varepsilon,\rho)\cup \mathbb{T}_\rho$ is a doubly
connected domain. One can easily see that $(D_1,D_2)\in
\mathcal{D}(E)$. Let $\rho_k(z)|dz|$ denote the extremal metric
for the corresponding module problem  for the domain $D_k$,
$k=1,2$; see \cite[Chapter II]{Jenkins}.

Since $D_2^*\subset D_2$, it follows that $\rho_2(z)|dz|$ is
admissible for the module problem for $D_2^*$ but it is not
extremal for this problem. Therefore, %
\begin{align} \label{11} %
m(D_2^*)&<\iint_{D_2} \rho_2^2(z)\,dA-\iint_{\cup_{k=1}^N
\Delta_k} \rho_2^2(z)\,dA
\\ %
{ }& =m(D_2)- \iint_{\cup_{k=1}^N \Delta_k} \rho_2^2(z)\,dA.
\nonumber
\end{align} %

Next, we consider a metric $\rho(z)|dz|$ on the domain $D_1^*$
defined by %
$$ 
\rho(z)=\begin{cases}  %
\rho_1(z)& {\mbox{if $z\in D_1$,}}\\
\rho_2(z)& {\mbox{if $z\in D_1^*\setminus D_1$.}} %
\end{cases} %
$$ 

We claim that $\rho(z)|dz|$ is admissible for the problem on the
reduced module for the domain $D_1^*$. To prove this, we consider
an analytic Jordan curve $\gamma\subset D_1^*$ separating
$0$ from $\partial D_1^*$. If $\gamma\subset \bar{D}_1$, then %
\begin{equation} \label{13} %
\int_\gamma \rho(z)\,|dz|=\int_\gamma \rho_1(z)\,|dz|\ge 1
\end{equation} %
since $\rho_1(z)|dz|$ is admissible for the problem on the reduced
module of $D_1$.

Suppose now that $\gamma\not\subset \bar{D}_1$. Then the circle
$\mathbb{T}_{\rho-\varepsilon}$ divides $\gamma$ into a finite
number of arcs. By $\tau_k$, $k=1,\ldots,m$, we denote those of
them, which lie in $D_2$. Let $\sigma_k\subset D_1^*$ denote the
closed arc of $T_{\rho-\varepsilon}$ joining the endpoints of
$\tau_k$, $k=1,\ldots,m$. It follows from Lemma~2, presented
below, that %
\begin{equation} \label{14} %
\int_{\tau_k}\rho(z)\,|dz|\ge \int_{\sigma_k}\rho(z)\,|dz|\ge
\int_{\sigma_k}\rho_1(z)\,|dz|.
\end{equation} %

Let $\tilde{\gamma}$ denote the curve obtained from $\gamma$ by
replacing the arcs $\tau_k$ with the arcs $\sigma_k$,
$k=1,\ldots,m$.  Then $\tilde{\gamma}$ is the curve in
$\bar{D}_1$, that is not Jordan in general, that separates $0$
from $\partial D_1$. Hence, %
\begin{equation} \label{15} %
\int_{\tilde{\gamma}}\rho_1(z)\,|dz|\ge 1.
\end{equation} %
Combining equations (\ref{14}) and (\ref{15}), we conclude that %
\begin{equation} \label{16} %
\int_{\gamma}\rho(z)\,|dz|\ge
\int_{\tilde{\gamma}}\rho_1(z)\,|dz|\ge 1.
\end{equation} %
Equations (\ref{13}) and (\ref{16}) imply that $\rho(z)|dz|$ is an
admissible metric for the problem on the reduced module in
$D_1^*$. Hence,

\begin{align} \label{17} %
m(D_1^*,0)&\le \lim_{\epsilon\to 0} \left\{\iint_{D_1^*\setminus
\bar{\mathbb{D}}_\epsilon} \rho^2(z)\,dA + \frac{1}{2\pi} \log
\epsilon\right\}
\\ %
{ }& =\lim_{\epsilon\to 0} \left\{\iint_{D_1\setminus
\bar{\mathbb{D}}_\epsilon} \rho^2(z)\,dA + \iint_{\cup_{k=1}^N
\Delta_k} \rho_2^2(z)\,dA+\frac{1}{2\pi} \log \epsilon\right\}
\nonumber \\ %
{ }&=m(D_1,0)+\iint_{\cup_{k=1}^N \Delta_k} \rho_2^2(z)\,dA.
\nonumber %
\end{align} %

Combining (\ref{11}) and ({\ref{17}), we obtain the following
inequalities:
\begin{align} 
\alpha_1^2m(D_1^*,0)&+\alpha_2^2m(D_2^*)<\alpha_1^2m(D_1,0)+\alpha_2^2m(D_2)\nonumber
\\  %
{ }&+(\alpha_2^2-\alpha_1^2)\iint_{\cup_{k=1}^N \Delta_k}
\rho_2^2(z)\,dA
\le \alpha_1^2m(D_1,0)+\alpha_2^2m(D_2). \nonumber %
\end{align} %
The latter inequalities contradict to the assumption that the pair
$(D_1^*,D_2^*)$ is extremal for Problem~M. This contradiction
shows that our assumption $R>r_2$ is wrong and therefore,
$L(E)\subset\bar{\mathbb{D}}_{r_2}$, as required.

The latter inclusion is proved in the case when $E$ consists of a
finite number of components. In the general case, we approximate
$E$ with a sequence of compact sets $E^j$, $j=1,2,\ldots$, such
that $E^{j+1}\subset E^j$ and such that $E^j$ is bounded by
appropriate level curves of Green's function
$g_{\mathbb{D}\setminus E}(z,0)$ of the domain
$\mathbb{D}\setminus E$ with pole at $0$. Then the required
inclusion $L(E)\subset \bar{\mathbb{D}}_{r_2}$  will follow from
the result already proved for sets $E$ with finite number of
connected components and from Carath\'{e}odory's convergence
theorem for simply connected and doubly connected domains.

The proof presented above, can be easily modified to show that if
$\alpha_1\ge \alpha_2$, then $L(E)\subset \bar{A}(r_1,1)$. %
 \hfill $\Box$

\begin{lemma} \label{Lemma 2}  %
Let $\sigma$ be an open arc on $\mathbb{T}$. Let $D_k$, $k=1,2$,
be a doubly connected domain having the circle $\mathbb{T}_r$ ,
$0<r<1$, as one of its boundary components. Suppose  that the
other boundary component of $D_1$, call it $\gamma_1$,  is such
that $l\subset \gamma_1$ and $D_1 \subset A(r,1)$. Suppose further
that  the other boundary component of $D_2$, call it $\gamma_2$,
is such that $l\subset \gamma_2$ and $A(r,1)\subset D_2$. Let
$\rho_k(z)|dz|$ denote  the extremal metric of the module problem
in $D_k$, $k=1,2$.
Then%
\begin{equation} \label{19} %
\rho_1(e^{i\theta})\le 1/2\pi\le \rho_2(e^{i\theta}) \quad
{\mbox{for all $e^{i\theta}\in \sigma$.}}
\end{equation} %
\end{lemma}  %

\noindent %
\emph{Proof.} Let us prove the first inequality. Let $f_1$ maps
$D_1$ conformally onto the annulus $A(r_1,1)$ such that
$f_1(\gamma_1)=\mathbb{T}$. Since the module of a doubly connected
domain increases under expansion, the following holds: $r<r_1<1$.
In terms of the mapping function, the extremal metric can be
expressed as follows (see \cite[Chapter II]{Jenkins}): %
$$  %
\rho_1(z)|dz|=\frac{1}{2\pi}\frac{|f'_1(z)|}{|f(z)|}\,|dz|.
$$ %
This shows that the desired result will follow if we prove that %
\begin{equation} \label{20}  %
{\mbox{meas}}(f_1(\sigma))\le {\mbox{meas}}(\sigma) %
\end{equation} %
 for every $\sigma\subset \mathbb{T}$, if the assumptions of
the lemma are satisfied.

To prove (\ref{20}), consider a family of curves
$\Gamma_1(\sigma)=\{\gamma\}$ consisting of all rectifiable arcs
$\gamma$ joining $\mathbb{T}_r$ and $\sigma$ inside the domain
$D_1$. Let $\Gamma(s,\sigma)$ denote a similar family of curves in
the annulus $A(s,1)$, $0<s<1$. Since the module increases under
expansion of the family of curves, we have %
\begin{equation}  \label{21} %
{\mbox{mod}}(\Gamma_1(\sigma))<{\mbox{mod}}(\Gamma(s,\sigma))
\quad {\mbox{if $D_1\not= A(r,1)$.}}
\end{equation} %

Now, suppose by contradiction that ${\mbox{meas}}(f_1(\sigma))>
{\mbox{meas}}(\sigma)$. Since the module of a family of curves is
conformally invariant and since a family of shorter curves has
bigger module than a corresponding family of longer curves, we
have the following: %
\begin{equation}  \label{22} %
{\mbox{mod}}(\Gamma_1(\sigma))={\mbox{mod}}(\Gamma(r_1,f_1(\sigma)))>{\mbox{mod}}(\Gamma(r,f_1(\sigma)))\ge
{\mbox{mod}}(\Gamma(s,\sigma)),
\end{equation} %
contradicting  (\ref{21}). Thus, the inequality (\ref{20}) is
proved and therefore the first inequality in (\ref{19}) holds. The
proof of the second inequality in (\ref{19}) follows the same
lines. %
\hfill $\Box$

\medskip

Notice that (\ref{19}) holds for doubly connected domains $D_k$
having the circle $\mathbb{T}_r$ of arbitrary small radius $0<r<1$
as one of its boundary components. Therefore, taking the limit as
$r\to 0$, we conclude that (\ref{19}) remains valid if $D_1$
and/or $D_2$ is a simply connected domain containing $0$. In the
latter case, $\rho_k(z)|dz|$ will denote the extremal metric for
the problem on the reduced module $m(D_k,0)$.

In the case of simply connected domains, (\ref{19}) also follows
from the Loewner's lemma or from Carleman's expansion principle
for the harmonic measure.

\medskip

\noindent%
\textbf{Remarks.} (1) As Gumenyuk and Roth showed  in \cite{GR},
Theorem~1 on the extremality of circularly slit disks cannot be
extended to domains of connectivity greater than two. This result
sounds similar to the observation made in \cite{S1} about a
possibility to extend Theorem~2 stated above in this note to the
case of several barriers $E_k$ contained in non-overlapping annuli
$\bar{A}(r_k^1,r_k^2)$, $0<r_1^1\le r_1^2<r_2^1\le
r_2^2<\cdots<r_n^1\le r_n^2<1$; i.e. such that $E_k\subset
\bar{A}(r_k^1,r_k^2)$. Precisely, if $(D_1^*,\ldots,D_n^*)$ is the
configuration of non-overlapping
domains maximizing the weighted sum of moduli %
$$
m(D_1,0)+\sum_{k=2}^n m(D_k),
$$
then the free boundary $L(E_1,\ldots,E_n)=\cup_{k=1}^n \partial
D_k^*\cap(\mathbb{D}\setminus \cup_{k=1}^n E_k)$ is not
necessarily contained in the union $\cup_{k=1}^n
\bar{A}(r_k^1,r_k^2)$ of these annuli. Here, $D_1\subset
\mathbb{D}\setminus \cup_{k=1}^n E_k$ is a simply connected domain
containing $0$ and $D_k$, $k=2,\ldots,n$, is a doubly connected
domain in $\mathbb{D}\setminus \cup_{k=1}^n E_k$ separating the
set $\cup_{j=1}^{k-1}E_j$ from the unit circle $\mathbb{T}$ and
the set $\cup_{j=k}^n E_j$.

\smallskip

\noindent %
(2) Our proof of Theorem~2 relies on the technique, which uses
weighted sums of moduli of free families of curves developed by
J.~A.~Jenkins and others. As well known, a module of a doubly
connected domain $D$ with complementary components  $E_0$ and
$E_1$ is the reciprocal of the capacity, when $D$ is considered as
the field of the condenser with plates $E_0$ and $E_1$. This
observation suggests that an approach utilizing properties of
capacities and potential functions of condensers also can be used
to study the squeezing problem. An example, demonstrating how this
approach works, was given in the second proof of Theorem~1.3 in
\cite{S1}. The proof of Theorem~2 in \cite{GR} is a nice
demonstration how this approach based on the potential theory can
be used in the context of the squeezing problem.


\bibliographystyle{amsplain}

\begin{thebibliography}{10}


\bibitem{DGZ} F. Deng, Q. Guan, and L. Znang, \textit{Some properties of squeezing functions on
bounded domains.} Pasific J. Math. \textbf{257} (2012), no. 2,
319--341.

\bibitem{Dubinin} 
V. N. Dubinin,   \emph{Condenser capacities and symmetrization in
geometric function theory.} Translated from the Russian by Nikolai
G. Kruzhilin. Springer, Basel, 2014. xii+344 pp.

\bibitem{GR} P. Gumenyuk and O. Roth, \textit{On the squeezing function for finitely connected planar domains.}     arXiv:2011.13734 [math.CV]
[v1] Fri, 27 Nov 2020. 

\bibitem{Jenkins} 
J. A. Jenkins, \emph{ Univalent functions and conformal mapping.}
Ergebnisse der Mathematik und ihrer Grenzgebiete. Neue Folge, Heft
18. Reihe: Moderne Funktionentheorie Springer-Verlag,
Berlin-G\"{o}ttingen-Heidelberg 1958 vi+169 pp.

\bibitem{Jenkins2} 
J. A. Jenkins, \emph{On the existence of certain general extremal
metrics.} Ann. of Math. (2) \textbf{66} (1957), 440--453. 


 \bibitem{Kuzmina} %
 G. V. Kuz'mina, \emph{Moduli of families of curves and quadratic differentials.}
 A translation of Trudy Mat. Inst. Steklov. \textbf{139} (1980).
 Proc. Steklov Inst. Math. 1982, no. 1, vii+231 pp. 

\bibitem{LSY} %
K. Liu, X. Sun, and S.-T.~ Yau, \textit{Canonical metrics on the moduli space of Riemann surfaces. I.}
J. Differential Geom., \textbf{68(3)}, (204), 571--637.

\bibitem{NTT} %
T. W. Ng, C. C. Tang, and J. Tsai, \textit{The squeezing function
on doubly connected domains  via the Loewner differential
equation.} Math. Ann.,
https://doi.org/10.1007/s00208-020-020--02046-w, 2020.


\bibitem{RW} %
 E.~Reich and S. E. Warschawski, \emph{On canonical conformal maps
of regions of arbitrary connectivity.}
 Pacific J. Math. \textbf{10} (1960), 965--985. 

\bibitem{S1}  %
A. Yu. Solynin, \emph{Geometric properties of
extremal partitions and estimates for the moduli of families of
curves in an annulus.} Zap. Nauchn. Sem. S.-Peterburg. Otdel. Mat.
Inst. Steklov. (POMI) \textbf{204} (1993), Anal. Teor. Chisel i
Teor. Funktsii. \textbf{11}, 93--114; translation in J. Math. Sci.
\textbf{79} (1996), no. 5, 1327--1340. %

 \bibitem{S2} %
 A. Yu. Solynin, \emph{Moduli and extremal metric problems.}
Algebra i Analiz \textbf{11} (1999), no. 1, 3--86; English
translation in St. Petersburg Math. J. \textbf{11} (2000), no. 1,
1--65. 

\end{thebibliography}

\end{document}